\documentclass[12pt, reqno]{amsart}
\usepackage{amsmath}
\usepackage{amssymb}
\usepackage{epsfig}
\usepackage{graphicx}
\usepackage{color}
\usepackage{fullpage}
\usepackage{comment}
\definecolor{shadecolor}{gray}{0.875}
\usepackage{amscd}
\usepackage{stmaryrd}

\numberwithin{equation}{section}

\input xy
\xyoption{all}

\calclayout
\allowdisplaybreaks[3]

\theoremstyle{plain}
\newtheorem{prop}{Proposition}[section]

\newtheorem{theo}[prop]{Theorem}

\newtheorem{lemm}[prop]{Lemma}

\theoremstyle{definition}
\newtheorem{defi}[prop]{Definition}

\newtheorem{ques}[prop]{Question}
\newtheorem{conj}[prop]{Conjecture}

\newtheorem{rema}[prop]{Remark}

\newtheorem{exam}[prop]{Example}

\def\Br{\mathrm{Br}}

\def\Pic{\mathrm{Pic}}

\def\Pic{\mathrm{Pic}}

\makeatother
\makeatletter

\author{Sho Tanimoto}
\address{Department of Mathematics, Faculty of Science, Kumamoto University, Kurokami 2-39-1 Kumamoto 860-8555 Japan}
\address{Priority Organization for Innovation and Excellence, Kumamoto University}
\email{stanimoto@kumamoto-u.ac.jp}

\title[Movable Bend and Break for sections]{Movable Bend and Break\\ for sections of del Pezzo fibrations}

\begin{document}
\date{\today}

\begin{abstract}
This is a report of the author's talk at Algebra Symposium 2019 held at Tohoku University.
We discuss some improvements of Mori's Bend and Break for sections of del Pezzo fibrations over $\mathbb P^1$.
\end{abstract}

\maketitle

\section{Introduction}

One of the main objects in diophantine geometry is a del Pezzo surface $S$ defined over a number field $k$, i.e., a smooth geometrically integral projective surface such that the anticanonical divisor $-K_S$ is ample. 
One can show that the degree $d = (-K_S)^2$ takes an integer value between $1$ and $9$, and we have the following classification of del Pezzo surfaces over an algebraic closure $\overline{k}$: let $\overline{S} = S \otimes_k \overline{k}$ and we have
\begin{itemize}
\item when $d = 9$, $\overline{S}$ is isomorphic to $\mathbb P^2$;
\item when $d = 8$,  $\overline{S}$ is isomorphic to $\mathbb P^1 \times \mathbb P^1$ or the blow up of the plane at one point;
\item when $d = 7$, $\overline{S}$ is isomorphic to the blow up of the plane at two distinct points;
\item when $d = 6$, $\overline{S}$ is isomorphic to the blow up of the plane at three non-colinear points;
\item when $d = 5$, $\overline{S}$ is isomorphic to the blow up of the plane at four points such that any three points are not colinear;
\item when $d = 4$, $\overline{S}$ is isomorphic to a complete intersection of two quadrics in $\mathbb P^4$;
\item when $d = 3$, $\overline{S}$ is isomorphic to a smooth cubic surface in $\mathbb P^3$;
\item when $d = 2$, $\overline{S}$ is isomorphic to a double cover of the plane ramified along a smooth quartic;
\item when $d = 1$, $\overline{S}$ is isomorphic to a double cover of a quadric cone ramified along a degree $6$ curve.
\end{itemize}
Moreover even when $1 \leq d \leq 4$, the surface $\overline{S}$ is also realized as the blow up of the plane at $9-d$ general points. Thus del Pezzo surfaces are geometrically rather well-understood objects.

However their arithmetic aspects remain mysterious and there are a lots of arithmetic questions one can ask for these surfaces. Among them the most important problems focus on the set $S(k)$ of rational points. Some of the natural questions on $S(k)$ are:
\begin{enumerate}
\item Existence of rational points;
\item Density of rational points;
\item Asymptotic formulas for the counting functions of rational points.
\end{enumerate}

It is a well-known fact that any smooth conic satisfies the Hasse principle, i.e., the existence of raitonal points is equivalent to the existence of local points at every place of $k$. However this feature fails for del Pezzo surfaces and so far all failures of the Hasse principle have been explained using the machinery of Brauer-Manin obstructions initiated in \cite{Manin}. Moreover Colliot-Th\'el\`ene conjecture predicts that the set of rational points is dense in the Brauer-Manin set and there are extensive studies on this subject. Readers interested in this story should consult \cite{Tony}.

The Zariski density of the set of rational points on a del Pezzo surface has been well studied too. Indeed if Colliot-Th\'el\`ene conjecture is correct then this implies that $S(k)$ is Zariski dense as soon as it is non-empty. This prediction suggests that the set $S(k)$ is always Zariski dense for any degree $1$ del Pezzo surface as the base locus of the anticanonical linear system always contains a unique rational point. However, this prediction is also out of reach to be proved at this moment.

The third question is about the asymptotic formula for the counting function of rational points of bounded height on a del Pezzo surface after removing the contribution from the exceptional set and Manin's Conjecture formulated in a series of papers \cite{BM90}, \cite{Peyre}, \cite{BT98}, \cite{Peyre03}, \cite{Pey17}, and \cite{LST18} predicts an explicit asymptotic formula in terms of geometric invariants of the underlying variety. There are extensive studies of this conjecture for del Pezzo surfaces, but so far there is no single smooth cubic surface defined over a number field which Manin's Conjecture has been established. We recommend interested readers to take a look at a survey paper on this subject \cite{Bro07}.

\

In this survey paper we would like to consider a del Pezzo surface over a different field, e.g.,  over $\mathbb C(t)$.
In this setting, Graber-Harris-Starr established the following vastly general theorem:
\begin{theo}{\cite{GHS03}}
Let $B$ be a smooth projective curve defined over an algebraically closed field $k$ of characteristic $0$. Let $X$ be a smooth projective geometrically rationally connected variety defined over the function field $k(B)$. Then $X$ admits a $k(B)$-rational point.
\end{theo}
Thus any del Pezzo surface defined over $\mathbb C(t)$ admits a rational point. Proving this theorem for del Pezzo surfaces is not so hard and one can appeal to the classification of minimal rational surfaces over non-closed fields. See \cite{Has09} for more details.
Weak approximation also has been studied in great details attested by \cite{HT06}, \cite{HT08}, \cite{Xu12}, \cite{Xu12b}, \cite{Tian}. The situation for del Pezzo surfaces seems to be almost complete, but there are still some open cases of del Pezzo surfaces of degree $1$ and $2$.

However, the situation on Manin's Conjecture in the settings of function fields is not satisfactory, and there seem to be few results in this direction compared to the situation over number fields. Over $\mathbb C(t)$, using Batyrev's heuristics, Manin's Conjecture can be interpreted as the following geometric problems:

Suppose we have a del Pezzo surface $X$ defined over $\mathbb C(t)$ and we fix an integral model $\pi : \mathcal X \to \mathbb P^1$. We denote the space of sections for $\pi$ by $\mathrm{Sec}(\mathcal X/\mathbb P^1)$. By valuative criterion there is a bijection between $\mathrm{Sec}(\mathcal X/\mathbb P^1)(\mathbb C)$ and $X(\mathbb C(t))$. Then we ask
\begin{enumerate}
\item What are the dimension and the number of components of bounded height for $\mathrm{Sec}(\mathcal X/\mathbb P^1)$?
\item Does the space $\mathrm{Sec}(\mathcal X/\mathbb P^1)$ enjoy some cohomological stablity when it is ordered by height?
\end{enumerate}

In \cite{LT19} and \cite{LTJAG}, Brian Lehmann and I have started a systematic study of Question (1) in the settings of trivial families of smooth Fano varieties using the geometry of the invariants appearing in Manin's Conjecture, which is developed in a series of papers \cite{HTT15}, \cite{LTT14}, \cite{HJ16}, \cite{LTDuke}, \cite{Sen17}, \cite{LST18}, and \cite{LTsurvey}. We obtained a satisfactory answer for the dimension of components of $\mathrm{Sec}(\mathcal X/\mathbb P^1)$ and classify irreducible components of $\mathrm{Sec}(\mathcal X/\mathbb P^1)$ for most Fano $3$-folds of Picard rank $1$.

In \cite{LTdP19} we take one step further and analyze sections of del Pezzo fibrations over $\mathbb P^1$. We obtain a satisfactory answer for a question about the dimension of $\mathrm{Sec}(\mathcal X/\mathbb P^1)$ (Theorem~\ref{theo:expecteddim}) and we establish Movable Bend and Break which is an improvement of Mori's Bend and Break (Theorem~\ref{theo:MBB}).

Mori invented a technique called Bend and Break lemma in \cite{Mori82}, and this shows that if we deform a rational curve while fixing two points, then it breaks into the union of rational curves. However, it is very difficult to control Bend and Break, i.e., a resulting curve may have more than two components and the corresponding point in the moduli space may not be a smooth point. Lehmann and I conjecture that for a Fano fibration, if the anticanonical height of section is sufficiently large, then it can break into the union of two free rational curves. We call this conjectural technique to be Movable Bend and Break and establish this for sections of del Pezzo fibrations over $\mathbb P^1$ in \cite{LTdP19}.

In this survey paper we discuss results in \cite{LTdP19} as well as the following applications:
\begin{enumerate}
\item Batyrev's conjecture for sections of del Pezzo fibrations;
\item Irreducibility of the space of sections for certain del Pezzo fibrations;
\item Stabilization of Abel-Jacobi maps;
\item Geometric Manin's Conjecture
\end{enumerate}

Here is the plan of the paper: In Section~\ref{sec:dim} we discuss the dimension of components of $\mathrm{Sec}(\mathcal X/\mathbb P^1)$. The main theorem is Theorem~\ref{theo:expecteddim} which claims that outside of a proper closed subset on $\mathcal X$, the dimension of the space of sections coincides with the expected dimension. In Section~\ref{sec:MBB}, we discuss our Movable Bend and Break lemma for sections of del Pezzo fibrations (Theorem~\ref{theo:MBB}) and demonstrate its proof in the case of surfaces. In Section~\ref{sec:applications} we discuss several applications of Movable Bend and Break as listed above.

\

\noindent
{\bf Acknowledgements:}
The author would like to thank Brian Lehmann for collaborations helping to shape his perspective on moduli of rational curves.
The author also would like to thank Brian for comments on an early draft of this paper.
The author would like to thank the organizers of Algebra Symposium 2019 for an opportunity to give a talk.
Sho Tanimoto is partially supported by Inamori Foundation, by JSPS KAKENHI Early-Career Scientists Grant numbers 19K14512, and by MEXT Japan, Leading Initiative for Excellent Young Researchers (LEADER).

\section{Expected Dimension}
\label{sec:dim}

In this paper we adopt the following definition:
\begin{defi}
An algebraic fiber space $\pi : \mathcal X \to \mathbb P^1$ is a del Pezzo fibration if
\begin{enumerate}
\item $\mathcal X$ is a smooth projective $3$-fold;
\item the relative anticanonical divisor $-K_{\mathcal X/\mathbb P^1} = -K_{\mathcal X} + \pi^*K_{\mathbb P^1}$ is relatively ample.
\end{enumerate}
\end{defi}
Corti showed the existence of a model with the ample anticanonical divisor with mild singularities in \cite{Cor96}.

\

A central object in this paper is the space of sections $\mathrm{Sec}(\mathcal X/\mathbb P^1)$. It is an open subscheme of the Hilbert scheme and it consists of countably many irreducible components.

For each $C \in \mathrm{Sec}(\mathcal X/\mathbb P^1)$, we define the height of $C$ by
\[
h(C):= -K_{\mathcal X/\mathbb P^1}.C.
\]
This height satisfies the Northcott property, i.e., the number of components of $\mathrm{Sec}(\mathcal X/\mathbb P^1)$ parametrizing sections of bounded height is finite.

Here are the main questions we ask for the space of sections:
\begin{enumerate}
\item for each component $M \subset \mathrm{Sec}(\mathcal X/\mathbb P^1)$, what is the dimension of $M$? Does it coincide with the expected dimension?
\item What is the number of components parametrizing sections of height $\leq d$?
\end{enumerate}
Understanding these questions is one of key assumptions of Batyrev's heuristics for Manin's Conjecture for trivial Fano families over finite fields, and we further apply this heuristics to del Pezzo fibrations in \cite{LTdP19}, leading to a conjectural solution to Geometric Manin's Conjecture formulated in \cite{LTdP19}.

\

In \cite{LTdP19}, we obtain a satisfactory answer for Question (1). To understand this result, we recall some deformation theory of rational sections:
we fix a del Pezzo fibration $\pi : \mathcal X \to \mathbb P^1$. We say a section $f : C \to \mathcal X$ is free if we have
\[
f^*T_{\mathcal X} = \mathcal O(a_1) \oplus \mathcal O(a_2) \oplus \mathcal O(a_3)
\]
with $0 \leq a_1 \leq a_2 \leq a_3$. If a section $C$ is free, then it is a smooth point of $\mathrm{Sec}(\mathcal X/\mathbb P^1)$ and the unique component containing $C$ has dimension equal to the expected dimension, i.e., 
\[
-K_{\mathcal X/\mathbb P^1}.C + 2.
\]

Let $M \subset \mathrm{Sec}(\mathcal X/\mathbb P^1)$ be a component and $p: \mathcal U \to M$ be the universal family. We say that $M$ is dominant if the evaluation map $s : \mathcal U \to \mathcal X$ is dominant. It is known that if a component $M$ is dominant, then a general member of $M$ is free. Hence one can conclude that $M$ has the expected dimension.

Keeping these facts in mind, now we state one of our main results from \cite{LTdP19}:

\begin{theo}{\cite[Theorem 1.1]{LTdP19}}
\label{theo:expecteddim}
Let $\pi : \mathcal X \to \mathbb P^1$ be a del Pezzo fibration.
Then there exists a proper closed subset $V \subsetneq \mathcal X$ such that any component of $\mathrm{Sec}(\mathcal X/\mathbb P^1)$ parametrizing a non-dominant family of sections will parametrize sections in $V$. In other words any section not contained in $V$ will deform to cover $\mathcal X$.
\end{theo}
Thus we can understand the dimension of each component inductively on dimension of the locus which sections sweep out.

\begin{rema}
In \cite{LT19}, we prove a similar statement for any smooth weak Fano variety replacing the space of sections by $\mathrm{Mor}(\mathbb P^1, X)$, and the main ingredient is the proper closedness of the exceptional set for Fujita invariants developed in \cite{LTT14}, \cite{HJ16}, \cite{LT19}. These results are based on the boundedness of singular Fano varieties proved by Birkar in \cite{birkar16} and \cite{birkar16b}.

In the above theorem, the set $V$ is also related to the Fujita invariant though the actual construction is more involved compared to trivial family cases.
\end{rema}

\section{Movable Bend and Break lemma}
\label{sec:MBB}

In this section we discuss Movable Bend and Break for sections of del Pezzo fibrations which is established in \cite{LTdP19}.
First of all we recall Mori's Bend and Break lemma from \cite{KM98}:

\begin{lemm}[Mori's Bend and Break]{\cite[Corollary 1.7]{KM98}}
Let $X$ be a projective variety and $f : \mathbb P^1 \to X$ be a stable map.
Suppose that we deform $f$ while fixing two distinct points on $X$ so that the image sweeps out a surface.
Then $f$ degenerates to a stable map $g: C \to X$ in the moduli space of stable maps such that
\begin{enumerate}
\item $C$ is a tree of rational curves;
\item $C$ consists of at least two non-contracted components;
\item $g(C)$ contains two points we fix.
\end{enumerate}

\end{lemm}

This lemma has vast applications to problems in algebraic geometry. Here are samples of applications:
\begin{itemize}
\item Cone theorem for smooth projective varieties (\cite{Mori82});
\item Rationally connectedness of smooth Fano varieties (\cite{KMM92} and \cite{Cam92});
\item Boundedness of smooth Fano varieties (\cite{KMM92});
\item Irreducibility of the space of rational curves on general Fano hypersurfaces (\cite{HRS04}, \cite{RY19}).
\end{itemize}

However, as mentioned in the introduction, there are certain difficulties in controlling Bend and Break:
\begin{itemize}
\item a resulting curve may have more than two components;
\item a resulting stable map may be a singular point of the moduli space of stable maps.
\end{itemize}
To overcome these issues we propose the following conjecture in \cite{LTdP19}:

\begin{conj}[Movable Bend and Break]{\cite[Conjecture 7.1]{LTdP19}}
\label{conj:MBB}
Let $\pi : \mathcal X \to \mathbb P^1$ be a Fano fibration, i.e, it is an algebraic fiber space such that the generic fiber is a smooth Fano variety and $\mathcal X$ is smooth. Then there exists a constant $Q(\mathcal X)$ with the following property: suppose that $C$ is a movable section such that $-K_{\mathcal X/\mathbb P^1}.C \geq Q(\mathcal X)$. Then $C$ deforms as a stable map to the union of two free curves.
\end{conj}
Note that a tree of free rational curves is a smooth point of the moduli space $\overline{M}_{0,0}(\mathcal X)$ of stable maps, so in particular there is the unique component containing it.
Here is the main theorem of \cite{LTdP19}:

\begin{theo}{\cite[Theorem 8.1]{LTdP19}}
\label{theo:MBB}
Conjecture~\ref{conj:MBB} holds when $\pi: \mathcal X \to \mathbb P^1$ is a del Pezzo fibration, i.e., the relative dimension is $2$ and $-K_{\mathcal X/\mathbb P^1}$ is relatively ample.

\end{theo}

Moreover we explicitly give a bound for $Q(\mathcal X)$. Let us describe this bound:
we define the minimum height of sections by
\[
\mathrm{neg}(\mathcal X/\mathbb P^1) = \min \{-K_{\mathcal X/\mathbb P^1}.C | \, C \in \mathrm{Sec}(\mathcal X/\mathbb P^1)\}.
\]
Note that this is well-defined due to the Northcott property of the height.
For an integer $d \in \mathbb Z$, we let $\mathrm{maxdef}(d)$ to be the maximum dimension of any component $M \subset \mathrm{Sec}(\mathcal X/\mathbb P^1)$ parametrizing sections of height $d$.
When there is no section of height $d$ we formally set $\mathrm{maxdef}(d) = - \infty$.
Then we define
\[
\mathrm{maxdef}(\mathcal X) = \max_{d < 0} \{\mathrm{maxdef}(d)\}.
\]
Here is our bound:
\begin{align*}
Q(\mathcal{X}) = \sup \{& 3, -2 \mathrm{neg}(\mathcal X/\mathbb P^1) -5, -\mathrm{neg}(\mathcal X/\mathbb P^1) + 3, \\ & 2\mathrm{maxdef}(\mathcal{X}) - 5\mathrm{neg}(\mathcal{X}/\mathbb{P}^{1}) - 5, \\
& 2\mathrm{maxdef}(\mathcal{X})-\mathrm{neg}(\mathcal{X}/\mathbb{P}^{1}) - 3, \\
& 2\mathrm{maxdef}(\mathcal{X}) + 2 + 2\sup\{0, - \mathrm{neg}(\mathcal{X}/\mathbb{P}^{1})\} \}.
\end{align*}
Furthermore, when $\mathrm{maxdef}(d)-d\leq 2$ holds for all $d < 0$ and there is no rational $-K_{\mathcal X_\eta}$-conic on the generic fiber $\mathcal X_\eta$ defined over $\mathbb C(\mathbb P^1)$, then Conjecture~\ref{conj:MBB} holds with $Q(\mathcal X) = 3$. (\cite[Lemma 8.1]{LTdP19})

\subsection{In dimension $2$}

Let us demonstrate a proof of Conjecture~\ref{conj:MBB} for surfaces:

\begin{prop}{\cite{LTdP19II}}
\label{prop:MBBforsurfaces}
Let $Y$ be a smooth projective surface with a morphism $\pi : Y \to \mathbb P^1$ such that a general fiber of $\pi$ is isomorphic to $\mathbb P^1$. Let $C$ be a section of $\pi$ such that
\[
-K_{Y/\mathbb P^1} \geq \max \{2, -\mathrm{neg}(Y/\mathbb P^1) + 1\}.
\]
Then one has
\[
C \sim C_0 + F
\]
where $C_0$ is a free section and $F$ is a general fiber of $\pi$.
\end{prop}

We will need the following lemma to prove Proposition~\ref{prop:MBBforsurfaces}:
\begin{lemm}{\cite{LTdP19II}}
\label{lemm:Lemma 5.3}
Let $F = \sum_{i} m_iE_i$ be a singular fiber of $\pi$ such that $m_1 = m_2 = 1$ and $E_i$'s are smooth rational curves. Let $Q = \sum_{i}a_i E_i \geq 0$ be an effective $\mathbb Q$-divisor such that we have
\[
Q.E_j = 
\begin{cases}
1 & \text{ if } j = 1\\
-1 & \text{ if } j = 2 \\
0 & \text{ otherwise }.
\end{cases}
\]
Then $-K_Y.Q > 0$.
\end{lemm}

This lemma can be proved using the induction on the number of components of $F$ and the MMP for surfaces.
Now let us explain a proof of Proposition~\ref{prop:MBBforsurfaces}:

\

{\bf A proof of Proposition~\ref{prop:MBBforsurfaces}}:
Since we have
\[
-K_{Y/\mathbb P^1}.C  = n \geq 2
\]
there exists a one-parameter family of deformations of $C$ passing through $n$ general points.
By Mori's Bend and Break we conclude that
\[
C \sim C_0 + mF + T,
\]
where $C_0$ is a section, $F$ a general fiber and $T$ is an effective divisor supported on singular fibers such that $T$ does not contain any full fiber.
Mori's Bend and Break actually shows that one can find a degeneration with at least two components going through general points (\cite[Lemma 4.1]{LTdP19}), and this means that we can assume that $m \geq 1$.

Let $F_0$ be a singular fiber and $T_0$ be the sum of terms in $T$ supported on $F_0$.
If $T_0$ is non-zero, then the negativity lemma implies that there exists some component of $F_0$ with non-vanishing intersection with $T_0$. Hence we conclude that $T_0$ satisfies the intersection property in Lemma~\ref{lemm:Lemma 5.3}, and this implies that $-K_Y.T_0 > 0$.
The upshot is that we have $b = -K_Y.T > 0$ unless $T = 0$.

Now our assumption on the height implies that $-K_Y.C_0 \geq 0$.
Thus $C_0$ can contain at most $n-2m-b +1$ general points. Each general fiber can contain at most $1$ general point. Thus we obtain
\[
n-2m-b +1 + m \geq n.
\]
This is only possible when $m = 1$ and $b = 0$.
Thus our assertion follows.
\qed

\

An idea of our proof in dimension $3$ is similar to the above proof. In dimension $3$, we separate analysis into two cases based on whether the normal bundle of a free section is balanced or not. When it is balanced, a similar but more complicated proof as above works fine. When it is not balanced, a free section will sweep out a surface after fixing an appropriate number of general points. Then we may reduce our analysis to the case of a surface.

\section{Applications}
\label{sec:applications}

In this section we discuss multiple applications of Movable Bend and Break lemma (Theorem~\ref{theo:MBB}).

\subsection{Batyrev's conjecture}

First of all we would like to introduce the following conjecture of Batyrev:

\begin{conj}[Batyrev's conjecture]
Let $X$ be a smooth projective weak Fano variety and $L$ be an ample divisor on $X$. Then there exists a polynomial $P(d)$ in $d$ such that the number of components of $\mathrm{Mor}(\mathbb P^1, X)$ parametrizing rational curves of $L$-degree $\leq d$ is at most $P(d)$.
\end{conj}

This conjecture should be contrasted with known exponential upper bounds for the number of irreducible components of Chow varieties/Hilbert schemes. See, e.g., \cite[I.3.28 Exercise]{kollar}.

In \cite{LT19}, Lehmann and I took the first step towards proving the above conjecture:

\begin{theo}{\cite[Theorem 1.4]{LT19}}
Let $X$ be a smooth projective uniruled variety and $L$ a big and nef divisor on $X$. Fix a positive integer $q$ and let $\overline{M} \subset \overline{M}_{0,0}(X)$ be the union of all components which contain a chain of free curves whose components have $L$-degree at most $q$.
Then there exists a polynomial $P(d)$ which is an upper bound for the number of components of $\overline{M}$ of $L$-degree at most $d$.
\end{theo}

Thus essentially this theorem means that Movable Bend and Break for rational curves will imply Batyrev's conjecture for dominant components. Combining the above theorem with Theorem~\ref{theo:MBB} ($+$ Proposition~\ref{prop:MBBforsurfaces}) we obtain:

\begin{theo}{\cite[Corollary 1.5]{LTdP19}}
Let $\pi : \mathcal X \to \mathbb P^1$ be a del Pezzo fibration.
Then there exists a polynomial $P(d)$ such that the number of components of $\mathrm{Sec}(\mathcal X/\mathbb P^1)$ parametrizing sections of height $\leq d$ is bounded by $P(d)$.
\end{theo}

\subsection{Irreducibility of the space of sections}

Movable Bend and Break (Theorem~\ref{theo:MBB}) can be used to prove the irreducibility of the space of sections of fixed height using induction on height:

\begin{exam}{\cite[Section 8.2]{LTdP19}}
\label{exam:fano}
Let $Y$ be a smooth Fano $3$-fold with $\Pic(Y) = \mathbb ZH$ and $-K_Y = 2H$. The degree $H^3$ takes an integer value between $1$ and $5$, and for each degree, there is exactly one deformation type of Fano $3$-folds. Here is the list of deformation types for $3 \leq H^3 \leq 5$:
\begin{itemize}
\item When $H^3 = 5$, $Y$ is a codimension $3$ linear section of $\mathrm{Gr}(2, 5) \subset \mathbb P^9$.
\item When $H^3 = 4$, $Y$ is a complete intersection of two quadrics in $\mathbb P^5$.
\item When $H^3 = 3$, $Y$ is a cubic smooth hypersurface in $\mathbb P^4$.
\end{itemize}

Assume that $d = H^3$ is in this range. Let $\beta: \mathcal X \to Y$ be the blow up of the base locus $Z$ of a general pencil of hyperplane sections. Then $\mathcal X$ comes with $\pi: \mathcal X \to \mathbb P^1$ induced by the pencil and $\pi$ is a del Pezzo fibration of degree $d$. We denote the exceptional divisor by $E \to Z$ and the space of sections of height $h$ by $M_h$.

\

\noindent
{\bf Claim}: $M_h$ is irreducible.

\

Indeed, first of all one can prove the following facts:
\begin{itemize}
\item
The minimum height is $-1$ and sections of the minimum height are fibers of $E \to Z$.
Note that $Z$ is an elliptic curve so there is no rational curve in $E$ other than fibers.
\item
The exceptional set in Theorem~\ref{theo:expecteddim} is $E$. Thus the only non-dominant family of sections are the family of sections of height $-1$.
\item
By \cite[Lemma 8.1]{LTdP19}, Movable Bend and Break works for sections of height $\geq 3$.
\item
One can check by hands that for $h = -1, 0, 1, 2$, $M_h$ is irreducible.
\end{itemize}

Assume that $h \geq 3$.
Suppose that $M_i$ is irreducible for any $i < h$.
Let $\overline{M}_h \subset \overline{M}_{0,0}(\mathcal X)$ be the Zariski closure of $M_h$ in the moduli space of stable maps $\overline{M}_{0,0}(\mathcal X)$.
By Movable Bend and Break, one can conclude that $C_0 + F \in \overline{M}_h$ where $C_0$ is a free section and $F$ is a free vertical curve. By \cite[Theorem 8.2]{LTdP19}, one can assume that $C_0$ is a section of height $h-2$ and $F$ is a vertical free conic.

Let $N$ be the space of vertical free conics -- this is irreducible because the monodromy on $N^1(\mathcal X_t)$ is the full Weyl group.
For any component $M \subset \overline{M}_{0,0}(\mathcal X)$, we denote, by $M^{(1)} \subset \overline{M}_{0,1}(\mathcal X)$, the family above $M$. We thus already proved that 
\[
C_0 + F \in \overline{M}_h \cap (M_{h-2}^{(1)} \times_{\mathcal X} N^{(1)}).
\]
Due to the maximal monodromy $M_{h-2}^{(1)} \times_{\mathcal X} N^{(1)}$ is also irreducible. On the other hand we have $C_0 + F \in \overline{M}_{0,0}(\mathcal X)^{\mathrm{sm}}$. Thus we conclude that
\[
M_{h-2}^{(1)} \times_{\mathcal X} N^{(1)} \subset \overline{M}_h
\]
It follows that $M_h$ is unique.

\end{exam}

\subsection{Stabilization of the Abel-Jacobi maps}

Let $\pi : \mathcal X \to \mathbb P^1$ be a del Pezzo fibration. Since $\mathcal X$ is rationally connected, the intermediate Jacobian
\[
\mathrm{IJ}(\mathcal X) = H^{2, 1}(\mathcal X)^\vee/H_3(\mathcal X, \mathbb Z)
\]
is a principally polarized abelian variety.
For any component $M \subset \mathrm{Sec}(\mathcal X)$ we have the Abel-Jacobi map
\[
\mathrm{AJ}_M : M \dashrightarrow \mathrm{IJ}(\mathcal X).
\]
When this map is dominant, we take a resolution $\beta: \widetilde{M} \to \overline{M}$ of a projective compactification $\overline{M}$ of $M$  such that $\mathrm{AJ}_M \circ \beta: \widetilde{M} \to \mathrm{IJ}(\mathcal X)$ is a morphism. We take the Stein factorization $Z_M \to \mathrm{IJ}(\mathcal X)$ of $\mathrm{AJ}_M \circ \beta$ and call it as the Stein factorization of $\mathrm{AJ}_M$.

Next theorem establishes the stabilization of the Abel-Jacobi maps when they are ordered by height.

\begin{theo}{\cite[Theorem 1.6]{LTdP19}}
\label{theo:stabilization}
Let $\pi : \mathcal X \to \mathbb P^1$ be a del Pezzo fibration of degree $\geq 3$.
Let
\[
\mathcal F = \{M \subset \mathrm{Sec}(\mathcal X/\mathbb P^1) \textnormal{ a component } | \textnormal{ $\mathrm{AJ}_M$ is dominant.}\}.
\]
Then the set
\[
\{Z_M \to \mathrm{IJ}(\mathcal X) | M \in \mathcal F \}/\textnormal{up to iso}
\]
is finite.
\end{theo}
A key lemma to prove the above theorem is the following:
\begin{lemm}{\cite[Proposition 10.1]{LTdP19}}
\label{lemm:stab}
In the settings of Theorem~\ref{theo:stabilization}, let $M \in \mathcal F$.
Let $M' \subset \mathrm{Sec}(\mathcal X/\mathbb P^1)$ be a component parametrizing smoothings of $C+F$ where $C \in M$ is a free section and $F$ is a free vertical conic or cubic.
Then $M' \in \mathcal F$.

Moreover let $Z_{M}\to \mathrm{IJ}(\mathcal X)$ and $Z_{M'} \to \mathrm{IJ}(\mathcal X)$ be the Stein factorizations of $\mathrm{AJ}_M$ and $\mathrm{AJ}_{M'}$ respetively.
Then $Z_{M}\to \mathrm{IJ}(\mathcal X)$ factors through as $Z \to Z' \to \mathrm{IJ}(\mathcal X)$.
\end{lemm}
Combining the above lemma with Theorem~\ref{theo:MBB}, one can deduce Theorem~\ref{theo:stabilization}. Let us illustrate this proof for Example~\ref{exam:fano}:

\begin{exam}{\cite[Example 10.4]{LTdP19}}
Let $Y$ be a smooth Fano $3$-fold with $\Pic(Y) = \mathbb ZH$, $-K_Y = 2H$ and $H^3 = 5$.
We define $\beta: \mathcal X \to Y$ as the blow up of the base locus of a general pencil of hyperplane sections on $Y$. Then $\mathcal X$ comes with a del Pezzo fibration $\pi: \mathcal X \to \mathbb P^1$ of degree $5$. Let $Z$ be the base locus of our pencil -- this is an elliptic curve. Since the intermediate Jacobian of $Y$ is trivial, we conclude that 
\[
\mathrm{IJ}(\mathcal X) = \mathrm{Jac}(Z).
\]
Let $M_h$ be the space of sections of height $h\geq 0$ -- we proved that this is irreducible.
Then the Abel Jacobi map $\mathrm{AJ}_{M_h}: M_h \dashrightarrow \mathrm{Jac}(Z)$ is described as follows:
For a general $C \in M_h$ its pushforward $\beta_*C$ is a rational curve of degree $h+1$ meeting $Z$ at $h$ distinct points.
Then the Abel-Jacobi map is described as
\[
M_h \dashrightarrow \mathrm{Sym}^h(Z) \rightarrow \mathrm{Jac}(Z), C \mapsto \beta_*C \cap Z \mapsto [\beta_*C \cap Z ] - [\beta_*C_0 \cap Z].
\]
Using this description one may prove that $\mathrm{AJ}_{M_i}$ is a MRC fibration for $i = 1, 2$.
Thus Lemma~\ref{lemm:stab} implies that $\mathrm{AJ}_{M_i}$ has connected fibers for any $i \geq 1$. 

It is natural to speculate the following questions:

\begin{ques}
Is $\mathrm{AJ}_{M_i}$ a MRC fibration for any $i \geq 1$?
\end{ques}

\begin{ques}
Can we conduct a similar analysis for $H^3 = 3, 4$?
\end{ques}

\end{exam}

\subsection{Geometric Manin's Conjecture}

Finally we discuss a conjectural solution to Geometric Manin's Conjecture from \cite{LTdP19}.
Let $\pi : \mathcal X \to \mathbb P^1$ be a del Pezzo fibration.
Due to smoothness of $\mathcal X$ every section intersects with a component of a $\pi$-vertical divisor with multiplicity $1$. Conversely Weak Approximation conjecture predicts that every possible intersection pattern of this type can be realized by some section. We call this intersection pattern as an intersection profile and we let $\Gamma_{\mathcal X}$ be the set of all possible intersection profiles, i.e., the set parametrizing the ways of choosing one component of multiplicity one for each fiber. Note that this is a finite set.

Let $N_1(\mathcal X)$ be the space of $\mathbb R$-cycles of dimension $1$ modulo numerical equivalence and $N_1(\mathcal X)_{\mathbb Z}\subset N_1(\mathcal X)$ be the lattice generated by integral cycles. Let $\mathrm{Nef}_1(\mathcal X) \subset N_1(\mathcal X)$ be the nef cone of curves.
Let $\lambda \in \Gamma_{\mathcal X}$ and we denote, by $N_\lambda \subset N_1(\mathcal X)$, the set of classes of a given intersection profile $\lambda$ -- this is an affine linear space of $N_1(\mathcal X)$.
Finally we denote the Brauer group of $\mathcal X$ by $\mathrm{Br}(\mathcal X)$.
The following conjecture is a key to our conjectural solution to Geometric Manin's Conjecture:

\begin{conj}{\cite[Conjecture 9.3]{LTdP19}}
\label{conj:irreducibility}
For each intersection profile $\lambda \in \Gamma_{\mathcal X}$. we let $\mathrm{Nef}_\lambda = \mathrm{Nef}_1(\mathcal X)\cap N_\lambda$. Then there is some translate $\mathcal T$ of $\mathrm{Nef}_\lambda$ in $N_\lambda$ such that every integral class in $\mathcal T$ is represented by exactly $|\mathrm{Br}(\mathcal X)|$ dominant families of sections.

\end{conj}

We now formulate Geometric Manin's Conjecture discussed in \cite{LT19} and \cite{LTdP19}. We assume that a general fiber of $\pi : \mathcal X \to \mathbb P^1$ is not isomorphic to $\mathbb P^2$ nor $\mathbb P^1\times \mathbb P^1$. First let us define the notion of Manin components:

\begin{defi}
We say a component $M\subset \mathrm{Sec}(\mathcal X/\mathbb P^1)$ a Manin component if for the universal family $p : \mathcal U \to M$, the evaluation map $s : \mathcal U \to \mathcal X$ does not factor through a proper subvariety $Y \subset \mathcal X$ such that the generic curve $Y_\eta$ is geometrically integral and has the anticanonical degree $\leq 2$. 
We let $\mathrm{Manin}_h$ be the set of Manin components parametrizing sections of height $h$.
\end{defi}

Using this we define the counting function which counts the number of Manin components of bounded height:
\begin{defi}
Fix a real number $q > 1$. For any positive integer $d$ we define the counting function
\[
N(\mathcal X, -K_{\mathcal X/\mathbb P^1}, q, d) = \sum_{h = 1}^d \sum_{M \in \mathrm{Manin}_h} q^{\dim M}.
\]
This function is closely related to the counting function over function fields of curves defined over $\mathbb F_q$.
\end{defi}

Geometric Manin's Conjecture predicts an asymptotic formula for $N(\mathcal X, -K_{\mathcal X/\mathbb P^1}, q, d)$ as $d \rightarrow \infty$. To demonstrate this formula let us define a few invariants of $\mathcal X$ which will be used to formulate the leading constant of the asymptotic formula:

\begin{defi}
We fix the Lebesgue measure $\mu$ on $N_1(\mathcal X_\eta)$ normalized so that the fundamental domain for $N_1(\mathcal X_\eta)_{\mathbb Z}$ has volume $1$. We define the alpha constant of $\mathcal X_\eta$ by
\[
\alpha(\mathcal X_\eta, -K_{\mathcal X/\mathbb P^1}) := \dim N_1(\mathcal X_\eta) \mu (\mathrm{Nef}(\mathcal X_\eta) \cap \{\gamma \in N_1(\mathcal X_\eta) | -K_{\mathcal X_\eta}.\gamma \leq 1\}).
\]
We also define the following invariant:
\[
\tau_{\mathcal X} = |\Gamma_{\mathcal X}|\cdot [N_1(\mathcal X)_{\mathbb Z}\cap N_1(\mathcal X_\eta): N_1(\mathcal X_\eta)_{\mathbb Z}]
\]
\end{defi}

Now let us state Geometric Manin's Conjecture for del Pezzo fibrations:

\begin{theo}{\cite[Theorem 9.10]{LTdP19}}
Let $\pi: \mathcal{X} \to \mathbb{P}^{1}$ be a del Pezzo fibration of degree $\geq 2$. 
Assume that Conjecture \ref{conj:irreducibility} holds for every intersection profile $\lambda\in \Gamma_{\mathcal X}$.  Then 
\begin{equation*}
N(\mathcal{X},-K_{\mathcal{X}/\mathbb{P}^{1}},q,d) \mathrel{\mathop{\sim}_{\mathrm{d \to \infty}}}  \frac{ \alpha (X_\eta, -K_{\mathcal{X}/\mathbb{P}^{1}})\tau_X|\Br(\mathcal{X})|}{1-q^{-1}}  q^{d} d^{\rho(\mathcal{X}_{\eta})-1}.
\end{equation*}
\end{theo}

The proof of this theorem is inspired by Batyrev's heuristics for Manin's conjecture over finite fields for trivial Fano families and uses Movable Bend and Break in an essential way.

\nocite{*}
\bibliographystyle{alpha}
\bibliography{AGSymp}

\end{document}